\documentclass[10pt]{article}
\usepackage{stmaryrd}
\usepackage{amsfonts}
\usepackage{amsmath}
\usepackage{amssymb}
\usepackage{bbm}
\usepackage{bm}
\usepackage{empheq}
\usepackage{latexsym}
\usepackage{amscd}
\usepackage{bbm}
\usepackage{url}
\usepackage{comment}
\usepackage{multirow}
\usepackage{xcolor}
\usepackage{dsfont}
\usepackage[left=3.2cm,right=3.2cm, top=2.5cm,bottom=2.5cm,bindingoffset=0cm]{geometry}

\newenvironment{Proof}{\par\noindent{\bf Proof of Theorem}} {\hfill$\scriptstyle\blacksquare$}

\newenvironment{ProofG}{\par\noindent{\bf Proof.}} {\hfill$\scriptstyle\blacksquare$}

\newcounter{theorem}[section]
\newtheorem{theorem}{Theorem}

\newtheorem{cor}[theorem]{Corollary}

\newtheorem{prop}[theorem]{Proposition}

\newtheorem{Remark}[theorem]{Remark}

\newcommand{\R}{\mathbb{R}}

\newcommand*{\Ind}[1]{\left\llbracket #1 \right\rrbracket}

\newcommand{\hash}{\scriptscriptstyle\#}

\title{Empirical forms of the Petty projection inequality}
\author{Grigoris Paouris \and Peter Pivovarov \and Kateryna Tatarko}

\date{\today}

\begin{document}

\maketitle

\begin{center}
    {\it  Dedicated to the memory of Clinton Myers
    Petty, 1923--2021.}
\end{center}

{\let\thefootnote\relax}  
{\let\thefootnote\relax}

\begin{abstract}
    The Petty projection inequality is a fundamental affine isoperimetric principle for convex sets. It has shaped several directions of research in convex geometry which forged new connections between projection bodies, centroid bodies, and mixed volume inequalities.  We establish several different empirical forms of the Petty projection inequality by re-examining these key relationships from a stochastic perspective.  In particular, we derive sharp extremal inequalities for several multiple-entry functionals of random convex sets, including mixed projection bodies and mixed volumes.
    
\end{abstract}

\section{Introduction}

\subsection{The Petty projection inequality and its reach}

Affine isoperimetric inequalities concern functionals on classes of
sets in which ellipsoids play an extremal role.  Typically such
inequalities involve convex bodies, taken modulo affine (or linear)
transformations, and are strictly stronger than their Euclidean
counterparts. The standard isoperimetric inequality can be derived
from several different affine strengthenings. Such affine inequalities have come to form an integral part of convex geometry and have
been extensively investigated within Brunn-Minkowski theory; see the expository survey \cite{Lutwak_survey} and books \cite{Gardner_book, Schneider_book} for foundational work on this subject.

A fundamental example is the Petty projection inequality. Recall that
the projection body $\Pi(K)$ of a convex body $K$ in $\mathbb{R}^n$ is
defined as follows: given a direction $\theta$ on the sphere
$S^{n-1}$, the support function of $\Pi(K) $ is the volume of
orthogonal projection of $K$ onto $\theta^{\perp}$ (see \S2 for
precise definitions). We write $ \Pi^{\circ}(K) $ for the polar of the
projection body. Petty's inequality states that among all convex
bodies of the same volume, ellipsoids maximize the volume of the polar
projection body. Formally, it can be stated as
\begin{equation}
\label{Petty}
| \Pi^{\circ} (K) | \leq | \Pi^{\circ} (K^*) |
\end{equation}  
where $K^*=rB_2^n$ is the centered Euclidean ball with radius $r$ chosen to
satisfy $|K^*|=|rB_2^n|$. The polar projection operator $\Pi^{\circ}$
satisfies $\Pi^{\circ} (TK)=T \Pi^{\circ}(K)$ for any
volume-preserving affine transformation $T$, which explains the use of
`affine' in this context.

Projection bodies are an important class of convex bodies in geometry
and functional analysis \cite{Bolker_1969,  BourLind_1988, Koldobsky_book, Schneider_Weil_Zonoids}. The volume of $\Pi^{\circ}(K)$ is
related to the surface area $S(K)$ via
$$\omega_n^{1/n}| \Pi^{\circ} (K) |^{-\frac{1}{n}} \leq S(K),$$ where
$\omega_n=|B_2^n|$; the latter follows directly from Cauchy's formula
and H\"{o}lder's inequality (see \cite[Remark 10.9.1]{Schneider_book}).  Thus Petty's
inequality implies the classical isoperimetric inequality for convex
sets. Up to normalization, the surface area $S(K)$ is one of the
quermassintegrals of $K$, while the quantity $|\Pi^{\circ}(K)|^{-1/n}$
is an {\it affine quermassintegral} of $K$.  Alexandrov's inequalities
state that among convex bodies of a given volume, all
quermassintegrals are minimized on balls (see \cite[\S 7.4]{Schneider_book}) . In a recent
breakthrough \cite{MY}, E. Milman and Yehudayoff proved that all affine
quermassintegrals are minimized on ellipsoids, verifying a
long-standing conjecture of Lutwak. This result establishes a
family of affine inequalities that interpolate between the Petty
projection inequality and the fundamental Blaschke-Santal\'{o}
inequality for the volume of the polar body of $K$.  The latter is
equalivalent to the affine isoperimetric inequality, see e.g., \cite{Lutwak_survey}.

Petty originally built on work of Busemann concerning the expected volume of
random simplices in convex bodies, and established what is known as
the Busemann-Petty centroid inequality \cite{Petty_centroid_1961}. He connected the latter to
projection bodies \cite{Petty_Isoperimetric_1967} by using an inequality about mixed volumes, known as
Minkowski's first inequality (\cite[\S 7.2]{Schneider_book}), which asserts that
\begin{equation*}
  V_1(K,L)\geq |K|^{n-1}|L|.
\end{equation*}
This idea was further developed by Lutwak and plays an important
role in kindred inequalities (see \cite{Lutwak_survey}). 

Since Petty's seminal work in 1972, his inequality has been proven by a number of different methods, e.g.,  \cite{Lutwak_survey, Sch_Petty_1995, Makai_Martini_Petty, MY}.  Moreover, several generalizations of the
inequality have been established. In particular, Lutwak, Yang and
Zhang introduced $L_p$ and Orlicz versions of the projection body and
proved the corresponding Petty inequalities \cite{LYZ_Lp_2000, LYZ_Orlicz_2010}. In \cite{Haberl_Schuster, Berg_Schuster}, a generalization to Minkowski valuations was obtained (see also \cite{Ludwig_2002} for a characterization of the projection body operator).  Another generalization, established by Lutwak involves the notion of mixed projection
bodies. Let $ K_{1}, \cdots , K_{n-1}$ be convex sets in $\mathbb
R^{n}$. The support function of the mixed projection body $ \Pi
(K_{1}, \cdots , K_{n-1})$ in a direction $\theta$ is defined as the following mixed volume:
\begin{equation}
   h_{\Pi(K_1,\ldots,K_{n-1})}(\theta) = nV(K_1,\ldots,K_{n-1},[0,\theta]),
\end{equation}
where $[0,\theta]$ is the line segment joining the origin and $\theta$.  Lutwak established several inequalities for mixed projection bodies,
one of which gives Petty's projection inequality as a special case \cite[Theorem 3.8]{Lutwak_1985_Mixed_projection}; namely,
\begin{equation}\label{eq::Lutwak_mixed}
   \lvert \Pi^{\circ}(K_1,\ldots,K_{n-1}) \rvert \leq \lvert \Pi^{\circ}(K_1^*,\ldots,K_{n-1}^*) \rvert.
\end{equation}
Recent active investigation around the notion of the projection body with respect to other measures and generalizations appear  in \cite{LRZ_measure_extension, LPRY}.

In this note we establish empirical versions of the Petty projection
inequality and its generalizations for mixed projection bodies. The
study of empirical versions of affine isoperimetric inequalities for centroid bodies and their $L_p$-analogues
was initiated by the first two authors in \cite{PaoPiv_2012} and further
developed in \cite{CEFPP_2015}. A number of inequalities in Brunn-Minkowski
theory have been shown to have stronger empirical forms \cite{PP_2017}, but Petty's
projection inequality has eluded our previous efforts. Inspired by recent results of E. Milman and Yehudayoff \cite{MY}, and also by the approach of Campi and Gronchi in \cite{Campi_Gronchi}, our work here is intended to fill this gap.

\subsection{Empirical mixed projection body inequalities}

Our main results concern randomly generated sets, obtained as linear images of a
compact, convex set $C\subseteq \mathbb{R}^m$ under an $n\times m$ random
matrix $\bm{X}$.  Namely, we will consider sets of the form
\begin{equation*}
  \bm{X} C = \left\{c_1X_1+\ldots+c_N X_N:(c_j)\in C\right\},
\end{equation*}
where $X_1,\ldots, X_m$ are independent random vectors distributed
according to densities of continuous probability distributions 
on $\mathbb{R}^n$.  We will write ${\bm X}^{\hash}$ for the $n\times
m$ random matrix that has independent columns distributed according to
$f^{*}$, the symmetric decreasing rearrangement of $f$ (see~\S~\ref{subsection:rearrangements}).

More generally, it will be convenient to work with matrices $\bm{X}$ whose column vectors are grouped into blocks.  Assume that $\{X_{ij}\}$ is a collection of independent random vectors such that $X_{ij}$ is distributed according to $f_{ij}$, $i=1,\ldots,n$, $j=1,\ldots,m_i$, where $m_i\geq 1$.  For $\ell=1,\ldots,n$, we write $m(\ell)=m_1+\cdots+m_{\ell}$,
and form $\bm{X} = [\bm{X}_1\ldots \bm{X}_{\ell}]$ with $n\times m_i$
blocks $\bm{X}_i=[X_{i1}\ldots X_{im_i}]$. We adopt a similar convention for ${\bm X}^{\hash}$, which consists of $n\times m_i$ blocks $\bm{X}_i^{\hash}=[X_{i1}^*\cdots X_{im_i}^*]$, where $X_{ij}^*$ are independent and distributed according to $f_{ij}^*$. For ease of reference, we summarize this notation in Table~\ref{table::notation}.
\begin{table}[h!]
    \centering
    \begin{tabular}{ |c|c|c| } 
\hline
\multirow{2}{6em}{$n\times m$ matrix with $\ell$ blocks } & \multirow{2}{6em}{$n\times m_i$ block,  $\ 1 \leq i \leq \ell$} & \multirow{2}{7em}{columns with densities} \\
&  & \\
\hline
$\bm{X} = [\bm{X}_1\ldots \bm{X}_{\ell}]$ & $\bm{X}_i=[X_{i1} \ldots X_{im_i}]$ & $X_{ij}\sim f_{ij}$\\
$\bm{X}^{\hash} = [\bm{X}_1^{\hash}\ldots \bm{X}_{\ell}^{\hash}]$ & $\bm{X}_i^{\hash}=[X_{i1}^*\ldots X^*_{im_i}]$ & $X^*_{ij}\sim f^*_{ij}$\\
\hline
\end{tabular}
    \caption{Random matrices with independent columns}
    \label{table::notation}
\end{table}

With this notation, our first main result concerns mixed projection bodies of random sets generated by $\bm{X}$ and $\bm{X}^{\hash}$.

\begin{theorem}\label{Petty_main}
  Let $C_1,\ldots,C_{n-1}$ be compact convex sets such that $\mathop{\rm dim}(C_i)=m_i$ for $i=1,\ldots,n-1$ and let $m=m_1+\ldots+m_{n-1}$.
  Let $\bm{X}$ and $\bm{X}^{\hash}$ be $n\times m$ random matrices with $\ell=n-1$ in Table~\ref{table::notation}. Then for any radial
  measure $\nu$ with a decreasing density,
  \begin{equation*}
  \mathbb{E} \nu \left(\Pi^{\circ}(\bm{X}C_1,\ldots,\bm{X}C_{n-1})\right) \leq \mathbb{E} \nu \left(\Pi^{\circ}(\bm{X}^{\hash}C_1,\ldots, \bm{X}^{\hash}C_{n-1})\right).
\end{equation*}
\end{theorem}

A special case of central importance concerns the classical projection body operator. Taking $\ell =1$ and writing $m=m(1)=m_1$ and $C=C_1=\ldots=C_n$, we have the following consequence.

\begin{theorem}
\label{thm:emp_Petty}
  Let $C$ be a compact convex set in $\mathbb{R}^m$. Let $\bm{X}$ and
  $\bm{X}^{\hash}$ be $n\times m$ random matrices with independent
  columns distributed according to $f$ and $f^{*}$, respectively. Then
  for any radial measure $\nu$ with a decreasing density,
  \begin{equation*}
    \mathbb E\nu(\Pi^{\circ}(\bm{X}C))\leq
    \mathbb E\nu(\Pi^{\circ}(\bm{X}^{\hash}C)).
  \end{equation*}
\end{theorem}

Theorem \ref{thm:emp_Petty}  extends the Petty projection inequality \eqref{Petty} in various
ways. Indeed, let $ K$ be a
convex body in $\mathbb R^{n}$ and let $ X_1, \ldots , X_m$ be
independent random vectors drawn uniformly from $K$. We denote their
convex hull by
\begin{equation*}
  [K]_m = \mathop{\rm conv}\{X_1,\ldots,X_m\}.
\end{equation*}
In matrix notation, we have $[K]_m= {\bm X}S_m $, where ${\bm
  X}=[X_1\dots X_m]$ and $S_m$ is the simplex $ S_m:= \mathop{\rm
  conv}\{ e_{1}, \dots , e_{m} \}$. Thus if $\nu$ is Lebesgue
measure, the above theorem states that
\begin{equation}
\label{emp_Petty}
 \mathbb E | \Pi^{\circ} ([K]_{m}) | \leq \mathbb E | \Pi^{\circ}
([K^*]_{m})|.
\end{equation}
Note that $\Pi^{\circ}([K^*]_m)$ is not a ball and the above statement
does not follow from Petty's inequality. However, when $m\rightarrow
\infty$, we get that $[K]_{m}\rightarrow K $, which implies
$\Pi^{\circ}([K]_m)\rightarrow \Pi^{\circ}(K)$, hence \eqref{Petty}
follows from \eqref{emp_Petty}.  The inequality $\nu\left( \Pi^{\circ}(K) \right) \leq \nu\left(\Pi^{\circ}(K^*) \right)$ that we get from Theorem~\ref{thm:emp_Petty} as $m \to \infty$ can also be directly obtained by adapting the proof of Petty's inequality in \cite[Section 8.2]{MY}.

More generally, let $K$ and $L$ be convex bodies in $\mathbb{R}^n$ and assume that the columns of ${\bm X}_1$ and ${\bm X}_2$ are distributed according to $\frac{1}{\lvert K \rvert}\mathds{1}_{K}$ and $\frac{1}{\lvert L\rvert}\mathds{1}_L$. For $p_1, p_2 \geq 1$, we define
 $[K]_{m_1}^{p_1}=\bm{X}_1 B_{p_1}^{m_1}$
and $[L]_{m_2}^{p_2}=\bm{X}_2 B_{p_2}^{m_2}$, we have
\begin{align}
    \label{eqn:rand_p_sum}
   \mathbb E \nu(\Pi^{\circ}([K]_{m_1}^{p_1}+_p[L]_{m_2}^{p_2}))\leq \mathbb E
   \nu(\Pi^{\circ}([K^*]_{m_1}^{p_1}+_p[L^*]_{m_2}^{p_2})),
\end{align}where $+_p$ denotes $L_p$-addition of sets $p\geq 1$ (see \S~\ref{sec::preliminaries}); in fact, we can accommodate more general $M$-addition and Orlicz addition operations (see \S~\ref{sec::preliminaries}). In a similar manner, Theorem \ref{Petty_main} implies 
\begin{equation}
\label{eq::ext_Lutwak}
 \mathbb E \nu( \Pi^{\circ} ([K]_{m_1}^{p_1}, \ldots, [K]_{m_{n-1}}^{p_{n-1}}) ) \leq \mathbb E \nu( \Pi^{\circ}
([K^*]_{m_1}^{p_1}, \dots, [K^*]_{m_{n-1}}^{p_{n-1}}))
\end{equation}
where we used the same notation as above. When $\nu$ is the Lebesgue measure, then \eqref{eq::ext_Lutwak} can be seen as a local version of \eqref{eq::Lutwak_mixed} for natural families of random convex sets associated to $K_1,\ldots,K_{n-1}$.

Further specializing to the case when $p_1=\ldots=p_{n-1}=\infty$, we get a corollary for the mixed projection body of centroid bodies. Recall that the centroid body $Z(L)$ of a convex body $L$ in $\R^n$ is defined by its support function via
\begin{align*}
    h_{Z(L)}(u) =\frac{1}{\lvert L \rvert}\int_L \lvert\langle x,u \rangle \rvert dx.
\end{align*} The Busemann-Petty centroid inequality mentioned above is a sharp extremal inequality for the volume of $Z(K)$, which heavily influenced affine isoperimetric principles \cite{Lutwak_survey} and the development of $L_p$-Brunn--Minkowski theory \cite{LYZ_Lp_2000,LYZ_Orlicz_2010}.

A stochastic notion of centroid bodies was developed in \cite{PaoPiv_2012}, which defines a random variant $Z_{m}(L)$ of $L$ as the body with support function 
\begin{align*}
    h_{Z_{m}(L)}(u)= \frac{1}{m}\sum_{i=1}^m \lvert \langle X_i,u \rangle \rvert,
\end{align*} where $X_1,\ldots,X_m$ are independent and identically distributed according to the normalized Lebesgue measure on $L$. Our next result concerns mixed projection bodies of independent empirical centroid bodies $Z_m(K_1),\ldots,Z_m(K_{n-1})$ whose support function is given by 

\begin{align*} 
    h_{\Pi(Z_m(K_1),\ldots,Z_m(K_{n-1}))}(y) &= nV(Z_m(K_1),\ldots,Z_m(K_{n-1}),[0,y]).   
\end{align*}

\begin{cor} \label{cor::Haddad}
Let $K_1, \dots, K_{n-1}$ be convex bodies in $\R^n$ and let $Z_m(K_1),\ldots,Z_m(K_{n-1})$ be independent empirical centroid bodies. Then 
\begin{equation*}
\mathbb E \nu(\Pi^{\circ}(Z_m(K_1),\ldots,Z_m(K_{n-1})))\leq
\mathbb E \nu(\Pi^{\circ}(Z_m(K_1^*),\ldots,Z_m(K_{n-1}^*))).
\end{equation*}
\end{cor}
We note that when $m \to \infty$,
\begin{equation} \label{eq::det}
   V(Z_m(K_1),\ldots,Z_m(K_{n-1}),[0,y]) \to \tilde{c}_n\int_{K_1}\cdots \int_{K_{n-1}}\lvert \mathop{\rm det}[x_1,\ldots x_{n-1},y]\rvert dx_1\ldots dx_{n-1} 
\end{equation}
where we use \cite[eq. (5.81)]{Schneider_book}. Haddad recently established a family of isoperimetric inequalities for a new class of convex bodies \cite{Haddad} that are defined using similar determinantal expressions as in \eqref{eq::det} and their $L_p$-generalizations. Our work shows that such bodies arise naturally as limiting cases of mixed projection bodies of random sets in Theorem \ref{Petty_main} when the $C_i$'s are chosen to be cubes.


\subsection{Empirical mixed volume inequalities}

We also present an alternate empirical version of Petty's projection inequality. This approach is inspired by the proof of the
inequality based on Busemann-Petty centroid inequality and 
Minkowski's first inequality \cite{Lutwak_survey, Schneider_book}. We will use an empirical approximant of centroid
bodies, defined as follows. For each convex body $L$ in
$\mathbb{R}^n$, we use the notation $[L]_{m_2}^{\infty}=\bm{X}_2 B_{\infty}^{m_2}$; this is nothing but the  Minkowski sum of $m_2$ random
segments $[-X_{2j},X_{2j}]$, where $X_{21},\ldots,X_{2m_2}$ are independent random
vectors sampled according to $\frac{1}{\lvert L \rvert}\mathds{1}_L$, i.e.,
\begin{equation*}
  [L]_{m_2}^{\infty} = \sum_{j=1}^{m_2}[-X_{2j}, X_{2j}].
\end{equation*}
Note that $\frac1m [L]_{m}^\infty = Z_{m}(L)$. Using this notation with $L=\Pi^{\circ}(K)$, we establish a sharp extremal inequality for the following quantity:
\begin{align*}
  &\mathbb E V_1([K]_{m_1},[\Pi^{\circ}(K)]_{m_2}^{\infty}) \\
    &=\frac{1}{\lvert K \rvert^{m_1} }\frac{1}{\lvert \Pi^{\circ}(K) \rvert^{m_2}}
    \int\limits_{K^{m_1}}\int\limits_{(\Pi^{\circ}(K))^{m_2}}V_1\left(\mathop{\rm conv}\{x_{11},\ldots,x_{1m_1}\}, \sum_{j=1}^{m_2}[-x_{2j},x_{2j}] \right)dx_{11}\ldots dx_{1m_1}
    dx_{21}\ldots dx_{2m_2};
\end{align*}
here we implicitly assume that the $m_1$ random vectors from $K$ and $m_2$ random vectors from $\Pi^{\circ}(K)$ are independent. With these notational conventions, we have the following theorem.

\begin{theorem} \label{Emp-Petty-2}
  Let $K$ be a convex body in $\mathbb{R}^n$ and $m_1,m_2\geq n$. Then
  \begin{equation*}
    \mathbb E V_1([K]_{m_1},[\Pi^{\circ}(K)]_{m_2}^{\infty}) \geq \mathbb E
    V_1([K^*]_{m_1},[\left(\Pi^{\circ}(K)\right)^*]_{m_2}^{\infty}).
  \end{equation*}
\end{theorem}

As we will explain at the end of \S~\ref{sec::LLN}, when $m_1,m_2\rightarrow \infty$, Theorem~\ref{Emp-Petty-2} also implies Petty's inequality \eqref{Petty}.

For the proof of Theorem \ref{Emp-Petty-2} we first need to
establish an empirical version of Minkowski's first inequality which
we believe is of independent interest. In fact, we establish a generalization of the latter, stated as follows.

\begin{theorem}
\label{Emp-Mixed}
  Let $C_1,\ldots,C_n$ be compact convex sets such that $C_i\subseteq
  \mathbb{R}^{m_i}$, $m_i\geq 1$, and set $m=m_1+\ldots+m_{n}$.
  Let $\bm{X}$ and $\bm{X}^{\hash}$ be $n\times m$ random matrices with $\ell=n$ in Table~\ref{table::notation}. Then
  \begin{equation*}
    \mathbb E V({\bm X}C_1,\ldots,{\bm X}C_n)\geq \mathbb E V({\bm
      X}^{\hash}C_1,\ldots,{\bm X}^{\hash}C_n).
  \end{equation*}
  Consequently, for any $k=1,\ldots,n-1$, 
  \begin{equation*}
    \mathbb E V({\bm X}C_1,\ldots,{\bm
      X}C_{k};B_2^n,n-k)\geq \mathbb E V({\bm
      X}^{\hash}C_1,\ldots,{\bm X}^{\hash}C_{k};B_2^n,n-k).
    \end{equation*}
\end{theorem}

Theorem \ref{Emp-Petty-2} follows directly from the latter theorem since 
\begin{align*}
\mathbb E &V_1([K]_{m_1},[\Pi^{\circ}(K)]_{m_2}^\infty)
    = \mathbb E V_1(\bm{X}_1S_{m_1},{\bm X}_2 B^{m_2}_\infty).
\end{align*}

In each of Theorems \ref{Petty_main}, \ref{Emp-Petty-2} and \ref{Emp-Mixed}, we have used a single matrix ${\bm X}$ with columns arranged in blocks  ${\bm X}_i$
according to Table \ref{table::notation}, and multiple bodies $C_1,\ldots,C_n$.
When the $C_i$ are all equal to a given compact convex set $C$ in $\mathbb{R}^{m_1}$, we have 
\begin{equation*}
V({\bm X}C_1,\ldots,{\bm X}C_n) = V({\bm X}_1 C,\ldots,{\bm X}_1 C) =\lvert {\bm X}_1C \rvert,
\end{equation*}
and only the first block ${\bm X}_1$ is involved in the expression; in particular, the block matrices
in the above mixed entry functional are dependent. When the $C_i$'s are compact convex sets placed in consecutive orthogonal subspaces $\mathbb{R}^{m_i}$, then the use of independent blocks ${\bm X}_i$ allows for distinct entries in \begin{equation*}
V({\bm X}C_1,\ldots,{\bm X}C_n)=V({\bm X}_1 C_1,\ldots,{\bm X}_n C_n)
\end{equation*} and all blocks ${\bm X}_i$ of ${\bm X}$ are used (and are independent).
The block notation for ${\bm X}$ also accommodates scenarios between these two extremes, where some of the $C_i$'s are repeated while others are taken in orthogonal subspaces.

\medskip

\noindent {\bf Acknowledgements.} We would like to thank Mark Rudelson and Franz Schuster for helpful discussions. The second and third authors are also grateful to ICERM for excellent working conditions, where they participated in the  program “Harmonic Analysis and Convexity”.

\medskip

\noindent {\bf Funding.} The first-named author was supported by NSF grants DMS 1800633, CCF 1900881, DMS 2405441, Simons Foundation Fellowship \#823432 and Simons Foundation Collaboration grant \#964286. The second-named author was supported by NSF grant DMS-2105468 and Simons Foundation grant \#635531. The third-named author was supported in part by NSERC
Grant no 2022-02961.


\section{Preliminaries} \label{sec::preliminaries}

\subsection{Convex bodies and mixed volumes}

We work in Euclidean space $\R^n$ and use $|\cdot|$ for the $n$-dimensional Lebesgue measure. The unit Euclidean ball in $\R^n$ is $B^n_2$, while the unit sphere is $S^{n-1}$. We use $P_{H}$ to denote the orthogonal projection onto a subspace $H$. We write  $u^\perp = \{x \in \R^n:\ \langle x,u \rangle = 0 \}$ for the $(n-1)$-dimensional subspace of $\R^n$ that is orthogonal to $u \in S^{n-1}$.

A convex body $K \subset\R^n$ is a compact, convex set with non-empty
interior. The set of all compact, convex sets is denoted by
$\mathcal{K}^n$. We say that $K$ is origin-symmetric if $K = -K$. We also say that $K$ is 1-unconditional if $K$ is symmetric with respect to reflections in the standard coordinate hyperplanes.

The support function of $K\in \mathcal{K}^n$ is defined by
$$
h_K(x) = \max\limits_{y \in K}\left< x, y \right>, \quad x \in \R^n.
$$
The polar body $K^\circ$ of an origin-symmetric convex body $K$ in $\mathbb{R}^n$ is defined as $K^\circ = \{x \in \R^n:\ h_K(x) \leq 1 \}$. 
The gauge function (or Minkowski functional) of an origin-symmetric convex body~$K$ is defined as $||x||_K = \inf \{ t \geq 0:\ x \in tK \}$. 
If $K$ contains the origin in its interior, then $||x||_K= h_{K^\circ}(x)$.

The Minkowski combination of $K_1, \dots, K_m \in \mathcal{K}^n$ is defined as
$$
\lambda_1 K_1 + \dots + \lambda_m K_m = \{\lambda_1 x_1 +\dots + \lambda_m x_m:\  x_i \in K_i\}
$$
where $\lambda_1, \dots, \lambda_m \geq 0.$
The Minkowski theorem on volume  of Minkowski combinations says that
$$
|\lambda_1 K_1 + \dots + \lambda_m K_m| = \sum\limits_{i_1, \dots, i_n = 1}^m \lambda_{i_1} \cdots \lambda_{i_n} V(K_{i_1}, \dots, K_{i_n}).
$$
The coefficient $V(K_{i_1}, \dots, K_{i_n})$ is the mixed volume of $K_{i_1}, \dots, K_{i_n}$; when the last body $K_{i_n}$ appears $\ell$ times, we write
$V(K_{i_1},\ldots,K_{i_{\ell}}; K_{i_n}, \ell)$, where $1\leq \ell \leq n$. For $K, L \in \mathcal{K}^n$ and $0 \leq i \leq n$, we write $V_i(K, L)$ to denote the mixed volume of $K$ repeated $n - i$ times and $L$ repeated $i$ times.

If $K_1, \dots, K_{n-1}$ are convex bodies in $\R^n$ and $u\in S^{n-1}$, then we write
$v(P_{u^\perp} K_1, \dots, P_{u^\perp} K_{n-1})$ for the
$(n-1)$-dimenisonal mixed volume of $P_{u^\perp} K_1, \dots,
P_{u^\perp} K_{n-1}$ in $u^\perp$. It is known (see
e.g. \cite[p. 302]{Schneider_book}) that for $u \in S^{n-1}$,
\begin{equation}\label{mixed_vol_proj}
v(P_{u^\perp} K_1, \dots, P_{u^\perp} K_{n-1}) = n V(K_1, \dots, K_{n-1}, [0, u])
\end{equation}
where $[0, u]$ denotes the line segment connecting the origin and $u$.

The projection body $\Pi K$ of a convex body $K$ is
defined as the origin-symmetric convex body such that $h_{\Pi K}(u) = |P_{u^\perp} K|$
for all $u\in S^{n-1}$. It follows from \eqref{mixed_vol_proj} that $h_{\Pi K} (u) = n V_1(K, [0,
  u])$ for $u \in S^{n-1}$. We will denote the polar projection body
$(\Pi K)^\circ$ by $\Pi^\circ K $.

More generally, the mixed projection body $\Pi(K_1, \dots, K_{n-1})$ of the convex
bodies $K_1, \dots, K_{n-1}$ is defined by
$$
h_{\Pi(K_1, \dots, K_{n-1})} (u) = v(P_{u^\perp} K_1, \dots, P_{u^\perp} K_{n-1}) = n V(K_1, \dots, K_{n-1}, [0, u])
$$ 
for all $u \in S^{n-1},$ where we used \eqref{mixed_vol_proj} in the last identity.
The polar of the mixed projection body $\Pi^\circ (K_1, \dots, K_{n-1})$ is defined as $(\Pi (K_1, \dots, K_{n-1}))^\circ$ so that
$$
||\theta||_{\Pi^\circ(K_1, \dots, K_{n-1})} =  n V(K_1, \dots, K_{n-1}, [0, \theta]).
$$For background on mixed projection bodies, see \cite{Lutwak_Mixed_Inequalities} and \cite{Schneider_book}.

\subsection{$L_p$ and $M$-addition operations}
\label{sec::Lp_addition}

We recall the notion of $L_p$-addition of convex bodies from $L_p$-Brunn--Minkowski theory, e.g. \cite{Firey,Lutwak93JDG,Lutwak96}. For
$K,L\in \mathcal{K}^n$ containing the origin and $p\geq 1$, we will
write $K+_{p}L$ for their $L_p$-sum, i.e.,
\begin{equation}
  \label{eqn:Msum}
  h^p_{K +_p L}(u)=h^p_{K}(u)+h^p_L(u),\quad u\in \mathbb{R}^n.
\end{equation}

A general framework for addition of sets, called $M$-addition,  was developed by Gardner, Hug and Weil \cite{GHW_2013}. Let $M$ be an arbitrary subset of $\mathbb{R}^m$ and define the $M$-combination  $\oplus_M (K_1, \dots, K_m)$ of arbitrary sets $K_1, \dots, K_m$ in $\mathbb{R}^n$ by

\begin{equation*}
\oplus_M (K_1, \dots, K_m) = \left\{ \sum_{i=1}^m a_i x_i : x_i \in K_i, \, (a_1, \dots, a_m) \in M \right\}.
\end{equation*}

If $M = \{(1, 1)\}$ and $K_1$, and $K_2$ are convex sets, then $K_1\oplus_M K_2 = K_1 + K_2$, i.e., $\oplus_M$ is the usual Minkowski addition. If $M = B_q^n, \ q \geq 1$ with $1/p + 1/q = 1$, and $K_1$ and $K_2$ are origin-symmetric convex bodies, then  $K_1 \oplus_M K_2 = K_1 +_p K_2$, i.e., $\oplus_M$ corresponds to $L_p$-addition as in \eqref{eqn:Msum}. More generally, let $\psi : [0, \infty)^2 \to [0, \infty)$ be convex, increasing in each argument, and  $\psi(0, 0) = 0$, $\psi(1, 0) = \psi(0, 1) = 1$. Let $K$ and $L$ be origin-symmetric convex bodies and let $M = B_\psi^\circ$, where $B_\psi = \{(t_1, t_2) \in [-1, 1]^2 : \psi(|t_1|, |t_2|) \leq 1 \}$. Then we define $K +_\psi L$ to be $K \oplus_M L$. In this way, $M$-addition encompasses previous notions of Orlicz addition, e.g. \cite{LYZ_Orlicz_2010}. 

It was shown in \cite[Section 6]{GHW_2013} that when $M$ is 1-unconditional and $K_1,\ldots, K_m$ are origin-symmetric and convex, then $\oplus_M (K_1, \dots, K_m)$ is origin-symmetric and convex. For our purposes, for such $M$ and $C_1,\ldots,C_m$, we have
\begin{equation}
\label{eqn:Madd}
  \oplus_M({\bf X}_1 C_1,\ldots,{\bf X}_n C_n)= [{\bf X}_1 \cdots {\bf X}_n]\oplus_{M}(C_1,\ldots,C_n);
\end{equation}see \cite[Sections 4 and 5]{PP_2017} for further background, details and references.

\subsection{Symmetrization of sets and functions}

\label{subsection:rearrangements}

Let $K \in \mathcal{K}^n$ and $u \in S^{n-1}$. We define $f_K:P_{u^{\perp}}K\rightarrow \R$ by
	\begin{equation*}
		\label{upperfunction}
		f_K(x)= \sup\{\lambda:x+\lambda u\in K\}
	\end{equation*}
	and $g_K:P_{u^{\perp}}K\rightarrow \R$ by
	\begin{equation*}
		\label{lowerfunction}
		g_K(x)=\inf\{\lambda:x+\lambda u\in K\}.
	\end{equation*}
	Notice that $-f_K$ and $g_K$ are convex functions. 
 
 The Steiner symmetral of a non-empty Borel set
	$A\subseteq\R^n$ of finite measure with respect to $u^{\perp}$, denoted here by
	$S_{u^{\perp}}A$, is constructed as follows: for each line
	$l$ orthogonal to $u^{\perp}$ such that $l\cap A$ is non-empty and measurable, the set
	$l\cap S_{u^{\perp}}A$ is a closed segment with midpoint on
	$u^{\perp}$ and length equal to the one-dimensional measure of $l\cap A$. 
 In particular, if $K$ is a convex body
	\begin{equation*}
		S_{u^{\perp}}K = \left\{x+\lambda u:x\in P_{u^{\perp}}K,
		-\dfrac{f_K(x)-g_K(x)}{2} \leq \lambda \leq
		\dfrac{f_K(x)-g_K(x)}{2}\right\}.
	\end{equation*}
	This shows that $S_{u^{\perp}}K$ is convex, since the function
	$f_K-g_K$ is concave. Moreover, $S_{u^{\perp}}K$ is symmetric with
	respect to $u^{\perp}$, it is closed, and by Fubini's theorem it
	has the same volume as $K$.

For a Borel set $A \subset \R^n$ with finite volume, the symmetric
rearrangement $A^*$ of $A$ is the open Euclidean ball centered at
the origin whose volume is equal to the volume of $A$. The symmetric
decreasing rearrangement of $\mathds{1}_A$ is defined as $\mathds{1}_A^* =
\mathds{1}_{A^*}$. It will be convenient to use the following bracket
notation for indicator functions:
\begin{equation}
  \label{eqn:indicator}
  \mathds{1}_A(x) = \Ind{x\in A}.
\end{equation}

Let $f: \R^n \rightarrow \R_+$ be an integrable function. Its
layer-cake representation is given by
\begin{equation}
  f(x) =\int_0^{\infty}\Ind{x\in\{f >t\}} dt.
\end{equation}
The symmetric decreasing rearrangement $f^*$ of $f$ is defined by
\begin{equation*}
    f^*(x) = \int\limits_0^\infty \Ind{x\in\{f>t\}^*} dt.
\end{equation*}
The function $f^*$ is radially-symmetric, radially decreasing and
equimeasurable with $f$, i.e. $\{f > t \}$ and $\{f^* > t \}$ have
the same volume for each $t>0.$ For integrable functions $f$, the Steiner symmetral $f^*(\cdot|u)$ of
	$f$ with respect to $u^{\perp}$ is defined as follows: 
\begin{align*}
f^*(x|u)=\int_{0}^{\infty}\Ind{x\in S_{u^{\perp}}\{f>t\}}
\end{align*}In other words, we obtain $f^*(\cdot|u)$ by rearranging $f$ along every line parallel to $u$. For more background on rearrangements, see \cite{LiebLoss_book}.

\section{Convexity and rearrangement inequalities}

\subsection{Shadow systems}

Rogers--Shephard \cite{RogShep_1958} and Shephard \cite{Shephard_1964} systematized the use of Steiner symmetrization as a means of proving geometric inequalities with their introduction of linear parameter systems and shadow systems, respectively.  A linear parameter system
is a family of sets
\begin{equation}
  \label{eqn:linear}
  K_{t}=\mathop{\rm conv}\{x_i+\alpha_i t u:i\in
  I\},
\end{equation}
where $\{\alpha_i\}_{i\in I}$ and $\{x_i\}_{i\in I}$ are bounded sets, and $I$ is an index set.  For a unit vector $u\in \mathbb{R}^n$ and a convex body $\mathcal{C}$
in $\mathbb{R}^{n+1}=\mathbb{R}^n\oplus \mathbb{R}$, a shadow system
is a family of sets of the form
\begin{equation*}
K_{t} = P_{t}\mathcal{C}, 
\end{equation*}
where $P_{t}:\mathbb{R}^{n+1}\rightarrow \mathbb{R}^n$ is the
projection parallel to $e_{n+1}- t u$. Setting 
\begin{equation*}
\mathcal{C}=
\mathop{\rm conv}\{x+\alpha(x) e_{n+1}:x\in K\subseteq e_{n+1}^{\perp}\}
\end{equation*}
where $\alpha$ is a bounded function on $K$, gives rise to the shadow system for $t\in[0,1]$,
\begin{equation*}
K_t=\mathop{\rm conv}\{x+t\alpha(x):x\in K\}.
\end{equation*}
The choice of $\alpha(x)=-g(P_{u^{\perp}}x)-f(P_{u^{\perp}}x)$ has the property that
$K_0=K$, while $K_{1/2}$ is the Steiner symmetral of $K$ about $u^{\perp}$, and $K_1$ is the reflection of $K$ about $u^{\perp}$. For background on linear parameter systems and shadow systems, we refer to \cite{Schneider_book,Campi_Gronchi}.

We will make essential use of the following fundamental theorem of Shephard 

\begin{theorem}
  \label{thm:Shephard}
Let $K_t^1,\ldots,K_t^n$ be shadow systems in common direction
$u$. Then
\begin{equation*}
  [0,1]\ni t\mapsto V(K_t^1,\ldots,K_t^n)
\end{equation*} is convex.

\end{theorem}


\subsection{Analytic tools}\label{sec::analytic_tools}

A non-negative, non-identically zero function $f$ is called
log-concave if $\log f$ is concave on $\{f > 0\}$. We note that if $f$
is a convex function, then the function $\Ind{f(x)\leq 1}$ is
$\log$-concave. Also, $f$ is quasi-concave if for all $t$ the set $\{x: f(x) > t\}$ is convex, and $f$ is quasi-convex if for all $t$ the set $\{x: f(x) \leq t\}$ is convex.

The Pr\'ekopa--Leindler inequality states that for $0 < \lambda <1$ and
functions $f, g, h: \R^n \rightarrow \R_+$ such that for any $x, y \in
\R^n$
$$ h(\lambda x + (1-\lambda) y) \geq f(x)^\lambda g(y)^{1-\lambda},
$$
the following inequality holds
$$
\int_{\R^n} h \geq \left(\int_{\R^n} f \right)^{\lambda}  \left(\int_{\R^n} g \right)^{1- \lambda}.
$$ We will use the following consequence of the Pr\'ekopa--Leindler
inequality: if $f:\mathbb{R}^n\times \mathbb{R}^m \to \R_+$ is $\log$-concave, then 
\begin{equation*}
  g(x) = \int_{\mathbb{R}^m}f(x,y)dx
\end{equation*} is a $\log$-concave function on $\mathbb{R}^n$; see \cite{Kanter}.

We also make use of Christ's form \cite{Christ} of the Rogers--Brascamp--Lieb--Lutinger inequality, see the survey \cite{PP_2017} for the related inequalities, their applications and further references.


\begin{theorem}\label{GCC}
Let $f_1, \dots, f_m: \mathbb{R}^n \rightarrow \mathbb{R_+}$ be integrable functions and let $F:(\mathbb{R}^n)^m \rightarrow \mathbb{R}_+$. Suppose that $F$ satisfies the following condition:
for any $u \in S^{n-1}$ and for any $\omega = (\omega_1, \dots, \omega_m) \subset~u^\perp$, the function $F_{u, \omega}: \mathbb{R}^m \rightarrow \mathbb{R}_+$ defined by $F_{u, \omega}(t_1, \dots, t_m) = F(\omega_1 + t_1 u, \dots, \omega_m + t_m u)$ is even and quasi-concave.
Then
\begin{align*}
\int\limits_{\mathbb{R}^n} \cdots \int\limits_{\mathbb{R}^n} F(x_1,
\dots, x_m) \prod_{i=1}^{m} f_i(x_i) dx_1 \dots dx_m \leq
\int\limits_{\mathbb{R}^n} \cdots \int\limits_{\mathbb{R}^n} F(x_1,
\dots, x_m) \prod_{i=1}^{m} f^*_i(x_i) dx_1 \dots dx_m .
\end{align*}
\end{theorem}
When each $F_{u, \omega}$ is even and quasi-convex, then the inequality in Theorem~\ref{GCC} is reversed. 

For subsequent reference, we note that the theorem is proved by iterated Steiner symmetrization and the key step involves Steiner symmetrization as follows:
\begin{align*}
  &\int_{\mathbb{R}^n} \cdots \int_{\mathbb{R}^n} F(x_1,
  \dots, x_m) \prod_{i=1}^{m} f_i(x_i) dx_1 \dots dx_m  =\\
  &=\int_{(u^{\perp})^m}\int_{\mathbb{R}^m}F_{u,\omega}(t_1, \dots, t_m)\prod_{i=1}^m f_i(\omega_i+t_iu)
  dt_1\dots dt_m  d\omega_1\dots d\omega_m\\
  & \leq \int_{(u^{\perp})^m}\int_{\mathbb{R}^m}F_{u,\omega}(t_1, \dots, t_m)\prod_{i=1}^m f^*_i(\omega_i+t_iu|u)
    dt_1\dots dt_m  d\omega_1\dots d\omega_m\\
    &=    \int_{\mathbb{R}^n} \cdots \int_{\mathbb{R}^n} F(x_1,
  \dots, x_m) \prod_{i=1}^{m} f^*_i(x_i|u) dx_1 \dots dx_m,
\end{align*}where $f_i^*(\cdot|u)$ is the Steiner symmetral of $f_i$ in direction $u$.

\section{Minimizing the mixed volume of random convex sets}
\label{section:mixed}

The proof of Theorem \ref{Emp-Mixed} relies on Theorem \ref{thm:Shephard} about the convexity of mixed volumes of shadow systems along a common direction.  Here we show how this interfaces with the use of random linear operators. 

\begin{Proof} {\bf \ref{Emp-Mixed}.} We start by associating a shadow system to linear images of convex sets. To fix the notation, let $C$ be a compact convex set in $\mathbb{R}^m$ and let $x_1,\ldots,x_m\in \mathbb{R}^n$. We will attach a shadow system in direction $u\in S^{n-1}$ to the set
\begin{equation*}
    [x_1\ldots x_m]C=\left\{\sum_{i=1}^m c_ix_i:(c_i)\in C\right\}.
\end{equation*}
 Decompose 
$x_i$ as $x_i=\omega_i+t_i u$, where $\omega_i\in
u^{\perp}$ for $i=1,\ldots,m$. Let ${\bm \omega}=(\omega_1,\ldots,\omega_m)$.  For ${\bf t}=(t_1,\ldots,t_m)\in
\mathbb{R}^m$, we form the $n\times m$ matrix
\begin{equation*}
T_{{\bm \omega}}({\bf t})=[\omega_1+t_1\theta \dots
  \omega_m+t_m\theta].
\end{equation*}
Fix ${\bf t},{\bf t}'\in \mathbb{R}^m$ and $c=(c_j)\in C$. Then for
each $\lambda\in[0,1]$,
\begin{eqnarray*}
T_{{\bm \omega}}(\lambda {\bf t}+(1-\lambda) {\bf t}')c &=&\sum_{j=1}^m
c_{j}(\omega_{j}+(\lambda t_{j}+(1-\lambda)t_{j}')u)\\ & = &
{\sum_{j=1}^mc_j(\omega_j+t_{j}'u)} +
\lambda{\Bigl(\sum_{j=1}^m
  c_i(t_{j}-t_{j}')\Bigr)}u.
\end{eqnarray*}
For $c\in C$, we write $x_c={\sum_{j=1}^mc_j(\omega_j+t_{j}'u)}$ and
$\alpha_c= \sum_{j=1}^m c_j(t_j-t_j')$ so that
\begin{equation}
  \label{eqn:image_shadow}
  T_{{\bm \omega}}(\lambda {\bf t}+(1-\lambda) {\bf
    t}')C=\left\{x_c+\lambda \alpha_c u:c\in C\right\}.
\end{equation}
As a linear image of the convex set $C$, the latter is convex and
hence takes the form of a linear parameter system in
\eqref{eqn:linear}, which is indexed by $C$ and generated by the
bounded sets $\{x_c\}_{c\in C}$ and $\{\alpha_c\}_{c\in C}$.

Similarly, assume we have $n$ compact convex sets $C_1,\ldots, C_n$ and $x_{ij}=\omega_{ij}+t_{ij}u$,
where $\omega_{ij}\in u^{\perp}$. For $i=1,\ldots,n$, let ${\bm
  \omega}_i=(\omega_{i1},\ldots,\omega_{im})$. For $\underline{{\bf t}}=({\bf t}_1,\ldots,{\bf t}_n)\in
\mathbb{R}^{m_1}\times \dots \times \mathbb{R}^{m_n}$, we write ${\bf t}_i=(t_{i1},\ldots,t_{im_i})\in
\mathbb{R}^{m_i}$ and set
\begin{equation*}
T_{{\bm \omega}_i}^i({\bf t}_i) = [\omega_{i1}+t_{i1}\theta \dots
  \omega_{im_i}+t_{im_i}\theta].
\end{equation*}

We first consider the case when $C_1, \ldots, C_n$ are in mutually orthogonal subspaces, then $\bm{X} C_i = \bm{X}_i C_i$ and the quantity under consideration is 
\begin{align}\label{eq::Thm1.4_quantity}
  \mathbb{E} &V({\bm X}_1C_1,\ldots,{\bm
    X}_{n}C_{n}) \nonumber \\
  &=\int_{\mathbb{R}^n}\cdots\int_{\mathbb{R}^n}
  V([x_{11}\cdots x_{1m_1}]C_1, \dots, [x_{n1}\cdots x_{nm_{n}}]C_{n}) \prod_{i=1}^{n} \prod_{j=1}^{m_i}  f(x_{ij})\,  dx_{11}\cdots dx_{nm_n}.
\end{align}

As we will apply Theorem \ref{GCC}, it is sufficient to show that
\begin{equation*}
  \mathbb{R}^{m_1}\times \cdots \times\mathbb{R}^{m_n}\ni ({\bf t}_1,\ldots,{\bf t}_n)\mapsto
  V(T_{{\bm \omega}_1}^1({\bf t}_1)C_1,\ldots,T_{{\bm \omega}_n}^n({\bf t}_n)C_n)
\end{equation*}
is convex. We need only show that the function is convex on any line joining given points $\underline{{\bf t}}=({\bf t}_1,\ldots,{\bf
  t}_n)$ and $\underline{{\bf t}}'=({\bf t}_1',\ldots,{\bf t}_n')$ in 
$\mathbb{R}^{m_1}\times \cdots \times\mathbb{R}^{m_n}$, i.e., we need only to establish convexity of
\begin{equation*}
[0,1]\ni \lambda \mapsto V(T_{{\bm \omega}_1}^1(\lambda {\bf
  t}_1+(1-\lambda){\bf t}_1')C_1,\ldots,T_{{\bm \omega}_n}^n(\lambda {\bf t}_n+(1-\lambda){\bf t}_n')C_n).
\end{equation*}

By the discussion at the beginning of this section, each argument in the mixed volume is a shadow system in the common direction $u$. Therefore, we can apply Theorem \ref{thm:Shephard} to obtain the required convexity in $\lambda$.

In the case $C_i$ are not necessarily mutually orthogonal, then $\bm{X}C_1, \ldots \bm{X}C_n$ share some common columns.
The proof then applies verbatim but on a smaller product space. For example,  in the case when all $C_i$ are identical, the mixed volume reduces to the volume and one works with 
\begin{equation*}
\int_{\mathbb{R}^n}\cdots
\int_{\mathbb{R}^n} \lvert[x_{11}\cdots x_{1m_1}]C_1\rvert
\prod_{j=1}^{m_1}f(x_{1j})dx_{11}\ldots dx_{1m_1};
\end{equation*} here the product space involves only the first $m_1$ random vectors.
As Theorem \ref{thm:Shephard} concerns any $n$ shadow systems (whether or not some are identical), it remains applicable in this case.  

\end{Proof}

\section{An empirical Petty projection inequality for measures}
\label{Proof_Petty}

While Theorem \ref{thm:emp_Petty} follows directly from Theorem \ref{Petty_main}, we will first prove the former for simplicity of exposition.  

\begin{Proof}~{\bf \ref{thm:emp_Petty}.}
We will first assume that $\nu$ is a  rotationally invariant, log-concave measure on
$\mathbb{R}^n$ with density $\rho (x)$, i.e. $d\nu(x) = \rho(x)
dx$. For $u \in S^{n-1}$ and $y \in u^\perp$, let $\rho_y(s) = \rho(y + su)$ be the restriction of $\rho$ to
the line that passes through $y$ parallel to $u$. We note
that since $\nu$ is a rotationally invariant, log-concave measure then
$\rho_y(s)$ is even and log-concave for any fixed $y \in u^\perp$.

Fix $u\in S^{n-1}$ and $\omega_1,\ldots,\omega_m\in u^{\perp}$. As
in \S\ref{section:mixed}, we write $T_{\bm \omega}({\bf t})=[\omega_1+t_1u \dots \omega_m+t_mu]$.
Using the notation for indicator functions \eqref{eqn:indicator}, we have
\begin{eqnarray}\label{volume}
\nu(\Pi^\circ(T_{\bm \omega}({\bf t})C)) &=&\int\limits_{u^\perp} \int\limits_{\R}
\Ind{y+su\in \Pi^\circ T_{\bm \omega}({\bf t})C} \rho_y(s) ds dy\\
&=&\int\limits_{u^{\perp}}\int\limits_{\R}\Ind{nV_{n-1}(T_{\bm \omega}({\bf t})C,[0,y+su])\leq 1}\rho_y(s) ds dy.
\end{eqnarray}
Fixed $y\in u^{\perp}$ and set
$h({\bf t}, s) = nV_{n-1}(T_{\bm \omega}({\bf t})C,\, [0, y+su])$, and  define
\begin{equation*}
G_{y}({\bf t}) = \int\limits_{\R} \Ind{h({\bf t}, s)\leq 1}\, \rho_y(s) ds.
\end{equation*}

Note that $h$ is jointly convex in ${\bf t}$ and $s$. To see this, it is
sufficient to show that, given any two points $({\bf t}, s)$ and $({\bf t}', s')$
in $\mathbb{R}^m\times \mathbb{R}$, the restriction of $h$ to the line
segment $[0,1]\ni \lambda\mapsto \lambda({\bf t},s)+(1-\lambda)({\bf t}',s')$ is
convex.  Set
\begin{equation*}
f(\lambda) = h(\lambda {\bf t} + (1-\lambda){\bf t}',
  \lambda s + (1-\lambda) s'),
\end{equation*} 
and observe that
\begin{eqnarray*}
 f(\lambda) =
 nV_{n-1}(T_{\bm \omega}(\lambda ({\bf t} - {\bf t}') + {\bf t}')C,\, [0, y+ (\lambda (s -
   s') + s')u]).
\end{eqnarray*}
Each of the arguments are shadow systems in a common direction $u$ as
observed above. Thus the convexity of $f(\lambda)$, and hence that of $h({\bf t}, s)$,  follows from Theorem~\ref{thm:Shephard}.

Using the joint convexity of $h$ in ${\bf t}$ and $s$, we have that
$\Ind{h({\bf t}, s) \leq 1}$ is $\log$-concave. As $G_y$ is the marginal of
$\log$-concave function on $\mathbb{R}^m\times \mathbb{R}$, it is also
$\log$-concave by the Prekopa--Leinder inequality (see \S \ref{subsection:rearrangements}).  In particular,
$G_y$ is quasi-concave.

Next, we note that $h$ is an even function. Indeed, the sets
$\left[\omega_1 + t_1 u, \dots, \omega_m + t_m u\right]C$ and
$\left[\omega_1 - t_1 u, \dots, \omega_m - t_m u\right]C$ are
reflections of one another with respect to $u^\perp$, hence
\begin{eqnarray*}
  h(-{\bf t}, -s) &= & nV_{n-1}(T_{\bm \omega}(-{\bf t})C,\, [0, y-su])\\ &=&
  n V_{n-1}(R_u(\left[\omega_1 + t_1 u, \dots, \omega_m + t_m
    u\right]C),\, R_u[0, y+su]) \\ &= & n V_{n-1}(\left[\omega_1 + t_1 u,
    \dots, \omega_m + t_m u\right]C,\, [0, y+su])\\&=& h({\bf t},s),
\end{eqnarray*}
where $R_u$ denotes the reflection with respect to $u^\perp$. The
rotational invariance of the density $\rho$ implies that $\rho_y$ is
even in $s$. Thus, by a change of variables, we have that $G_y$ is
even
\begin{eqnarray*}
  G_y(- {\bf t}) =\int\limits_{\R} \Ind{h(-{\bf t},s)\leq 1} \rho_y(s) ds =
  \int\limits_{\R} \Ind{h({\bf t},-s)\leq 1} \rho_y(s) ds
  =\int\limits_{\R} \Ind{h({\bf t},s)\leq 1} \rho_y(-s) ds = G_y({\bf t}).
\end{eqnarray*}

Thus for every $u \in S^{n-1}$ and for every $\bm \omega$, the
function $G_y({\bf t})$ is quasi-concave for every $y \in u^\perp$. We write
\begin{eqnarray*}
  \mathbb{E} \left( \nu( \Pi^{\circ}({\bm X} C)) \right) &=&
  \int_{\mathbb{R}^n}\cdots\int_{\mathbb{R}^n}
  \nu(\Pi^{\circ}([x_1\cdots x_m]C)) \prod_{i=1}^m f(x_i) \, dx_1\cdots
  dx_m \\ & = & \int_{(u^{\perp})^m}\int_{\mathbb{R}^m}
  \nu(\Pi^{\circ}(T_{\bm \omega}({\bf t})C)
  \prod_{i=1}^m f(\omega_i+t_iu)\, d{\bf t} d\overline{y}\\ & = &
  \int_{(u^{\perp})^m}\int_{\mathbb{R}^m} \int_{u^{\perp}}G_y({\bf t})dy
  \prod_{i=1}^m f(\omega_i+t_iu) \, d{\bf t}d\overline{y} \\ &=
  &\int_{(u^{\perp})^m} \int_{u^{\perp}}
  \int_{\mathbb{R}^m} G_y({\bf t})\prod_{i=1}^m f(\omega_i+t_iu) d{\bf t} dy
  d\overline{y}. \\
\end{eqnarray*}
Applying the key step in the proof of Theorem \ref{GCC} from \S \ref{sec::analytic_tools}, we obtain:
\begin{equation*}
  \mathbb{E} \left( \nu( \Pi^{\circ}({\bm X} C)) \right) \leq
      \mathbb{E} \left( \nu( \Pi^{\circ}({\bm X}^u C)) \right),
\end{equation*}
where ${\bm X}^u$ has independent columns distributed according to the
Steiner symmetrals $f^*(\cdot|u)$. By iterated symmetrizations, we
obtain
\begin{equation*}
  \mathbb{E} \left( \nu( \Pi^{\circ}({\bm X} C)) \right) \leq
    \mathbb{E} \left( \nu( \Pi^{\circ}({\bm X}^{\hash} C)) \right).
\end{equation*}

As the proof shows, Theorem~\ref{GCC} does not require densities to be identical. Hence, Theorem~\ref{thm:emp_Petty} holds under this less restrictive assumption.

Up to this point, we dealt with a log-concave measure $\nu$. The above applies, in particular, to the case of the Lebesgue measure restricted to a centered Euclidean ball. Next, we can consider a measure $\nu$ such that
$$
d\nu(x) = \rho(x) dx  \qquad \text{with} \qquad \rho: [0, \infty) \to [0, \infty) \qquad \text{decreasing}.
$$
Using Fubini's theorem, we have
$$
\mathbb{E}  \nu(\Pi^{\circ}({\bm X} C)) = \mathbb{E} \int\limits_0^\infty |\Pi^{\circ}({\bm X} C) \cap \{\rho \geq t\}| dt = \int\limits_0^{\rho(0)} \mathbb{E}  |\Pi^{\circ}({\bm X} C) \cap R(t)B^n_2| dt
 $$
where $R(t)$ is a radius of an Euclidean ball $\{\rho \geq t\}$ which implies the result for any radial measure $\nu$ with a decreasing density $\rho$.
\end{Proof}

\begin{cor}Let $K$ and $L$ be convex bodies in $\R^n$.  Let $M\subseteq \mathbb{R}^2$ be $1$-unconditional and compact.
Then for $p_1, p_2 \geq 1$
\begin{equation*}
    \mathbb{E} \nu(\Pi^{\circ}([K]_{m_1}^{p_1} \oplus_M [L]_{m_2}^{p_2}))\leq 
        \mathbb{E} \nu(\Pi^{\circ}([K^*]_{m_1}^{p_1} \oplus_M [L^*]_{m_2}^{p_2})).
\end{equation*}
\end{cor}

If $M=B_q^2$, with $1/p+1/q=1$, then $\oplus_M$ coincides with $+_p$-addition as defined in \S \ref{sec::Lp_addition}. Thus the corollary immediately implies \eqref{eqn:rand_p_sum}.

\begin{ProofG} Let the columns of ${\bm X}_1$ and ${\bm X}_2$ be independent and distributed according to $\frac{1}{\lvert K \rvert}\mathds{1}_{K}$ and $\frac{1}{\lvert L\rvert}\mathds{1}_L$, respectively.  As in the introduction, $[K]_{m_1}^{p_1}={\bm X}_1 B_{p_1}^{m_1}$ and $[L]_{m_2}^{p_2}={\bm X}_2 B_{p_2}^{m_2}$. Assuming $B_{p_1}^{m_1}\subseteq \mathop{\rm span}\{e_1,\ldots,e_{m_1}\}$ and writing $\widehat{B_{p_2}^{m_2}}$ for the copy of $B_{p_2}^{m_2}$ lying in $\mathop{\rm span}\{e_{m_1+1},\ldots,e_{m_1+m_2}\}$, we have
\begin{align*}
[K]_{m_1}^{p_1} \oplus_M [L]_{m_2}^{p_2} & ={\bm X}_1 B_{p_1}^{m_1} \oplus_{M} {\bm X}_2 B_{p_2}^{m_2}\\
&=[{\bm X}_1 {\bm X}_2](B_{p_1}^{m_1}\oplus_M \widehat{B_{p_2}^{m_2}})\\&=
{\bm X}(B_{p_1}^{m_1}\oplus_M \widehat{B_{p_2}^{m_2}}).
\end{align*}
In this way the $M$-addition operation coincides with the image of the convex body $C=B_{p_1}^{m_1}\oplus_M \widehat{B_{p_2}^{m_2}}$ under ${\bm X}$.  Similarly, we have 
\begin{align*}
[K^*]_{m_1}^{p_1} \oplus_M [L^*]_{m_2}^{p_2} = {\bm X}^{\#}(B_{p_1}^{m_1}\oplus_M \widehat{B_{p_2}^{m_2}}).
\end{align*}
Thus Theorem \ref{thm:emp_Petty} applies directly.
\end{ProofG}

\begin{Remark}The latter corollary is a special case of the following inequality
\begin{equation*}
\mathbb E\nu(\Pi^{\circ}(\oplus_M (\bm{X}_1 C_1,\ldots,\bm{X}_n C_n) ))\leq 
\mathbb E\nu(\Pi^{\circ}(\oplus_M (\bm{X}^{\#}_1 C_1,\ldots,\bm{X}^{\#}_n C_n))),
\end{equation*}which follows directly from Theorem \ref{thm:emp_Petty} and \eqref{eqn:Madd} whenever $\oplus_{M}(C_1,\ldots,C_n)$ is compact and convex.
\end{Remark}

\section{Mixed projection bodies}

In this section, we prove Theorem \ref{Petty_main}.  As the argument is similar to the proof given in \S~\ref{Proof_Petty}, we simply outline the additional steps. 

\begin{Proof}~{\bf\ref{Petty_main}.} Let $\nu$ be a log-concave, rotationally invariant measure on
$\mathbb{R}^n$ with a density $\rho (x)$ as in \S~\ref{Proof_Petty}. As in \S~\ref{section:mixed}, we start by assuming that $C_1, \ldots, C_n$ lie in the mutually orthogonal subspaces.

Fix $u\in S^{n-1}$ and $\bm{\omega}_{i} = (\omega_{i1}, \dots, \omega_{im})$ such that $\omega_{ij} \in u^{\perp}$ for $i = 1, \dots, n-1$. Using notation in \S~\ref{section:mixed}, we write $T^i_{\bm \omega_{i}}({{\bf t}_i})=[\omega_{i1}+t_{i1}\theta \cdots
  \omega_{im}+t_{im}\theta].$ Thus, we have
\begin{eqnarray}\label{volume_mixed_proj}
&&\nu(\Pi^\circ(T^1_{\bm \omega_1}({\bf t}_1)C_1, \dots, T^{n-1}_{\bm \omega_{n-1}}({\bf t}_{n-1})C_{n-1})) \nonumber \\
&=&\int\limits_{u^\perp} \int\limits_{\R}
\Ind{y + su \in \Pi^\circ(T^1_{\bm \omega_1}({\bf t}_1)C_1, \dots, T^{n-1}_{\bm \omega_{n-1}}({\bf t}_{n-1})C_{n-1})} \rho_y(s) ds dy \nonumber\\
&=&\int\limits_{u^{\perp}}\int\limits_{\R}\Ind{nV(T^1_{\bm \omega_1}({\bf t}_1)C_1, \dots, T^{n-1}_{\bm \omega_{n-1}}({\bf t}_{n-1})C_{n-1},[0,y+su])\leq 1}\rho_y(s) ds dy.
\end{eqnarray}

For fixed $y\in u^\perp$, we define 
$$
h(\underline{{\bf t}}, s) = V(T^1_{\bm \omega_1}({\bf t}_1)C_1, \dots, T^{n-1}_{\bm \omega_{n-1}}({\bf t}_{n-1})C_{n-1},[0,y+su])$$
and 
$$
G_{y}(\underline{{\bf t}}) = \int\limits_{\R} \Ind{h(\underline{{\bf t}}, s)\leq 1}\, \rho_y(s) ds.
$$

Using the same reasoning as in Section~\ref{Proof_Petty}, $h(\underline{{\bf t}}, s)$ is jointly convex in  $\underline{{\bf t}}$ and $s$, and an even function. In particular, joint convexity implies that $\Ind{h(\underline{{\bf t}}, s) \leq 1}$ is $\log$-concave. Also, we have that $G_y$ is
even, i.e., $G_y(-\underline{{\bf t}}) =  G_y(\underline{{\bf t}}).$ Therefore, for every $u \in S^{n-1}$ and for every $\bm \omega_i$, the function $G_y(\underline{{\bf t}})$ is quasi-concave for every $y \in u^\perp$.

Repeating the same argument on $(\mathbb{R}^n)^{m_1} \times \dots \times (\mathbb{R}^n)^{m_{n-1}}$ as in \S\ref{Proof_Petty}, we get 
\begin{equation*}
  \mathbb{E} \left( \nu( \Pi^{\circ}({\bm X}_1 C_1, \dots, {\bm X}_{n-1} C_{n-1})) \right) \leq
      \mathbb{E} \left( \nu( \Pi^{\circ}({\bm X}^{u}_1 C_1, \dots, {\bm X}^{u}_{n-1} C_{n-1}))  \right),
\end{equation*}
and after iteration of the repeated symmetrization in suitable directions, we arrive at the following
\begin{equation*}
  \mathbb{E} \left( \nu( \Pi^{\circ}({\bm X}_1 C_1, \dots, {\bm X}_{n-1} C_{n-1})) \right) \leq
    \mathbb{E} \left( \nu( \Pi^{\circ}({\bm X}^{\hash}_1 C_1, \dots, {\bm X}^{\hash}_{n-1} C_{n-1})) \right).
\end{equation*}

Once again, when $C_1, \ldots, C_n$ are non-mutually orthogonal, the argument applies verbatim on the smaller product space. 

Finally, as above, the proof shows that densities need not be identical. 
\end{Proof}

\section{Laws of large numbers and convergence}
\label{sec::LLN}

In this section, we detail how one can obtain deterministic inequalities from our main stochastic inequalities. As in the introduction, the centroid body $Z(L)$ of a convex body $L$ in $\R^n$,is defined by
\begin{align*}
    h_{Z(L)}(u) =\frac{1}{\lvert L \rvert}\int_L \lvert\langle x,u \rangle \rvert dx.
\end{align*}
Similarly, the empirical centroid body $Z_{1,N}(L)$ of $L$ is given by  
\begin{align*}
    h_{Z_{1,N}(L)}(u)= \frac{1}{N}\sum_{i=1}^N \lvert \langle X_i,u \rangle \rvert,
\end{align*}where $X_1,\ldots,X_N$ are independent and identically distributed according to the normalized Lebesgue measure on $L$. Note that by the strong law of large numbers (as in e.g., \cite{PaoPiv_2012}), we have convergence a.s. in the Hausdorff metric:
\begin{align}\label{eq:SLN_centroid}
    Z_{1,N}(L)\rightarrow Z(L).
\end{align}

\begin{prop} \label{prop::LLN}Let $K$ be a convex body in $\mathbb{R}^n$ and $m_1,m_2\geq n$. Then
    \begin{align}
 \frac{1}{m_2}\mathbb E V_1([K]_{m_1},[\Pi^{\circ}(K)]_{m_2}^{\infty})\rightarrow V_1(K,Z(\Pi^{\circ}(K))). 
    \end{align}
\end{prop}

\begin{ProofG}
    We have convergence a.s. in the Hausdorff metric
    \begin{align*}
        [K]_{m_1} = \mathop{\rm conv}\{X_{11},\ldots,X_{1m_{1}}\}\rightarrow K;
    \end{align*}see, e.g., \cite{DL_rates_convergence} for a stronger quantitative result on the rate of convergence.
 Similarly, by \eqref{eq:SLN_centroid}
    \begin{align*}
        \frac{1}{m_2}[L]_{m_2}^{\infty}=\frac{1}{m_2}\sum_{j=1}^{m_2}[-X_{2j},X_{2j}] \rightarrow Z(L).
    \end{align*}
It follows that we have convergence a.s., as $m_1,m_2\rightarrow \infty$,
\begin{align*}
   \frac1{m_2} V_1([K]_{m_1},[\Pi^{\circ}(K)]_{m_2}^{\infty})\rightarrow V_1(K,Z(\Pi^{\circ}(K))).
\end{align*}
Note that 
\begin{align*}
    V_1([K]_{m_1},[\Pi^{\circ}(K)]_{m_2}^{\infty})&\leq V_1(K,[\Pi^{\circ}(K)]_{m_2}^{\infty})\\&=\frac{1}{n}\int_{S^{n-1}}h_{[\Pi^{\circ}(K)]_{m_2}^{\infty}}(u) d\sigma_K(u)\\
    &= \frac{1}{n}\frac{1}{m_2}\sum_{j=1}^{m_2}\int_{S^{n-1}}
        \lvert \langle u, X_{2j}\rangle \rvert d\sigma_K(u)\\
        & \leq \frac{1}{n} R(\Pi^{\circ}(K))R(K)S(K)   
\end{align*}
where $R(K)$ is the circumradius of $K$ and $S(K)$ is the surface area of $K$.
By the bounded convergence theorem, we have as $m_1,m_2\rightarrow \infty$,
\begin{align*}
   \frac1{m_2} \mathbb E V_1([K]_{m_1},[\Pi^{\circ}(K)]_{m_2}^{\infty})\rightarrow V_1(K,Z(\Pi^{\circ}(K))).
\end{align*}

\end{ProofG}

\medskip

Lastly, we show that Thereom~\ref{Emp-Petty-2} implies Petty's projection inequality~\eqref{Petty}. Applying Theorem~\ref{Emp-Petty-2} as $m_1,m_2\rightarrow \infty$ and using Proposition~\ref{prop::LLN}, we get
\begin{equation} \label{eq::mixed_vol}
    V_1(K, Z(\Pi^{\circ}(K))) \geq V_1(K^*, Z((\Pi^{\circ}(K))^*)).
\end{equation}
Now we appeal to the fact that 
\begin{align*}
    V_1(K,Z(\Pi^{\circ}(K)))=c_n,
\end{align*}where $c_n$ is a constant independent of $K$; see, e.g., \cite{Lutwak_survey}, and compute the right-hand side of~\eqref{eq::mixed_vol}. For this, we first observe that $(\Pi^\circ K)^*  = \frac{|\Pi^\circ(K)|^{\frac1n}}{|\Pi^\circ(K^*)|^{\frac1n}} \Pi^\circ(K^*) $ (see \cite[Corollary 9.1.4]{Gardner_book}). Then 
\begin{align*}
    V_1(K^*, Z((\Pi^{\circ}(K))^*))
    &= \frac{|\Pi^\circ(K)|^{\frac1n}}{|\Pi^\circ(K^*)|^{\frac1n}}  V_1(K^*, Z(\Pi^{\circ}(K^*))) = c_n \frac{|\Pi^\circ(K)|^{\frac1n}}{|\Pi^\circ(K^*)|^{\frac1n}}
\end{align*}
where we used homogeneity of mixed volumes and $Z(\phi K) = \phi Z(K)$ for $\phi \in GL_n$. Thus,
\begin{align*}
c_n &\geq \frac{|\Pi^\circ(K)|^{\frac1n}}{|\Pi^\circ(K^*)|^{\frac1n}} c_n
\end{align*}
which implies \eqref{Petty}.
\begin{Remark}
    The above derivation closely resembles the proof of the Petty projection inequality using  the first Minkowski inequality and Busemann centroid inequality which can be found in \cite{Lutwak_survey}. However, the stochastic inequality in Theorem~\ref{Emp-Petty-2} in effect combines the use of these two inequalities into one step.
\end{Remark}

\bibliography{sample}

\bibliographystyle{plain}

\bigskip

\bigskip

G.~Paouris, \textsc{Department of Mathematics, Texas A\&M University, College Station, Texas, 77843-3368, USA}\par\nopagebreak
\textsc{Department of Mathematics,
Princeton University, Fine Hall, 304 Washington Road, Princeton, NJ, 08540, USA}\par\nopagebreak
  \textit{Email address}: \texttt{grigoris@math.tamu.edu}

\medskip

P.~Pivovarov, \textsc{Mathematics Department, University of Missouri, Columbia, Missouri, 65211, USA}\par\nopagebreak
  \textit{Email address}: \texttt{pivovarovp@missouri.edu}

\medskip

K.~Tatarko, \textsc{Department of Pure Mathematics, University of Waterloo, Waterloo, Ontario, N2L 3G1, Canada}\par\nopagebreak
  \textit{Email address}: \texttt{ktatarko@uwaterloo.ca}

\end{document}